\documentclass[12pt,a4paper]{article}
\usepackage{amsthm,amssymb,amsmath,mathrsfs,graphicx}
\usepackage[pagebackref=false,colorlinks,linkcolor=black,citecolor=magenta]{hyperref}
\allowdisplaybreaks[1]

\newtheorem{theorem}{\bf Theorem}
\newtheorem{proposition}{\bf Proposition}
\newtheorem{lemma}{\bf Lemma}[section]

\newtheorem{corollary}[lemma]{\bf Corollary}
\newtheorem{remark}{\bf Remark}

\renewcommand{\proof}{\noindent{\it\textbf{Proof.}}\ \ }

\newcommand{\eqd}{$\hfill \blacksquare$}

\newcommand{\Sym}{{\rm Sym}}

\date{}

\title {Finite nilpotent groups coincide with their $2$-closures in all of their faithful permutation representations}

\author{Alireza Abdollahi\footnote{E-mail: a.abdollahi@math.ui.ac.ir} $^{a,b}$ 
 and Majid Arezoomand\footnote{Email: arezoomand@math.iut.ac.ir} $^a$\\ {\small\em $^a$ Department of Mathematics, University of Isfahan, Isfahan 81746-73441, Iran}\\ {\small\em $^b$ School of Mathematics,} \\ {\small\em Institute for Research in Fundamental Sciences (IPM),} \\ {\small \em P. O. Box 19395-5746, Tehran, Iran}   }

\begin{document}
\maketitle
\noindent\textbf{Abstract.}

Let $G$ be any group and $G$ be a subgroup of $\Sym(\Omega)$ for some set $\Omega$.
The \textit{$2$-closure} of $G$ on $\Omega$,
 denoted by $G^{(2),\Omega}$, is by definition,
 \[\{\theta\in\Sym(\Omega)\mid \forall \alpha,\beta\in\Omega, \exists g\in G, \alpha^\theta=\alpha^g,\beta^\theta=\beta^g\}.\]

The group $G$ is called  {\it $2$-closed} on $\Omega$ if $G=G^{(2),\Omega}$. We say that
 a  group $G$ is a {\it $2$-closed group} if $H=H^{(2),\Omega}$ for any set $\Omega$ such that $G\cong H\leq\Sym(\Omega)$.
 Here we show that the center of any finite $2$-closed group is cyclic and
 a finite nilpotent  group is $2$-closed if and only if it is cyclic or a direct product of a generalized quaternion group with a cyclic group of odd order.


\section{Introduction and results}
Let $\Omega$ be a set and $G$ be a group with $G\leq\Sym(\Omega)$. Then $G$ acts naturally on
$\Omega\times\Omega$ by
$(\alpha_1,\alpha_2)^g=(\alpha_1^g,\alpha_2^g)$, where $g\in G$ and $\alpha_1,\alpha_2\in\Omega$.
The {\it $2$-closure} of $G$ on $\Omega$,
denoted by $G^{(2),\Omega}$, is
 \[\{\theta\in\Sym(\Omega)\mid \forall \alpha,\beta\in\Omega, \exists g\in G, \alpha^\theta=\alpha^g,\beta^\theta=\beta^g\}.\]
 $G^{(2),\Omega}$ contains $G$ and is the  largest subgroup
of $\Sym(\Omega)$ whose orbits on $\Omega\times\Omega$ are the same orbits of $G$ \cite[Definition 5.3 and Theorem 5.4]{Wielandt}.
The notion of $2$-closure as a tool in the study of permutation groups was introduced
by I. Schur \cite{Wielandt}. The study of $2$-closures of permutation groups has been initiated by Wielandt \cite{Wielandt} in 1969,
to present a
unified treatment of finite
and infinite permutation groups, based on  invariant relations and invariant functions.  For further studies and applications see
\cite{Faradzev, Liebeck0, Liebeck1, Onan, Ponomarenko1, Ponomarenko2,Praeger1, Praeger2, Xu}.
Since $G\leq G^{(2),\Omega}$, it is an interesting question that how far the $2$-closure of a permutation group $G$ is from $G$?
In 2012, K. M. Monks in his Ph.D. thesis proved that every finite cyclic group is $2$-closed in all of its faithful permutation representations \cite[Theorem 12]{Monks}.
Derek Holt in his answer to a question \cite{math}  proposed by  the second author in {\sf Mathoverflow} introduced a class of abstract groups called $2$-closed groups  consisting  of   all groups which are $2$-closed in all of their faithful permutation representations. That is, an abstract group $G$ is  {\it $2$-closed} if  $H=H^{(2),\Omega}$ for any  set $\Omega$ with $G\cong H\leq\Sym(\Omega)$. So all finite cyclic groups by  Monk's result are $2$-closed.  Holt \cite{math} mentioned that generalized quaternion groups $Q_{2^n}$ of order $2^n$, $n\geq 3$,  and the quasicyclic $p$-group $\mathbb{Z}_{p^\infty}$ are also $2$-closed and the group $Q_{2^n} \times C_2$ is not $2$-closed, where $C_m$ is a finite cyclic group of order $m$.  

 Our main results are the following.



\begin{theorem}\label{cycliccenter}
The center of every finite $2$-closed group is cyclic.
\end{theorem}

We determine all finite nilpotent $2$-closed groups.

\begin{theorem}\label{nilp} A finite nilpotent  group is $2$-closed if and only if it is cyclic or a direct product of a generalized quaternion group with a cyclic group of odd order.
\end{theorem}

 Our notations are standard and are mainly taken from \cite{Dixon}, but for
the reader's convenience  we recall some of them as  follows:
\begin{itemize}
\item[] $\Sym(\Omega)$: The symmetric group on the set $\Omega$.
\item[] $Q_{2^n}$: The generalized quaternion group of order $2^n$, $n\geq 3$.
\item[] $C_n$: The finite cyclic group of order $n$.
\item[] $\alpha^g$: The action of $g$ on $\alpha$.
\item[] $G_\alpha$: The point stabilizer of $\alpha$ in $G$.
\item[] $\alpha^G$: The orbit of $\alpha$ under $G$.
\item[] $Z(G)$: The center of $G$.
\item[] $Z_2(G)$: The second center of $G$.
\item[] $C_G(H)$: The centralizer of the  subgroup $H$ of a group $G$.
\item[] $N_G(H)$: The normalizer of the subgroup $H$ of a group $G$.
\item[] $H_G$: The core of the  subgroup $H$ of $G$, that is the intersection of all $G$-conjugates of $H$.
\end{itemize}

\section{Proof of Theorem 1}
Recall that a group $G$ is a {\it $2$-closed group} if $H=H^{(2),\Omega}$ for any set $\Omega$ such that $G\cong H\leq\Sym(\Omega)$. First, let us prove or recall
some useful results.
\begin{lemma}{\rm (\cite[Exercise 5.29]{Wielandt})}\label{commute}
Let $G\leq\Sym(\Omega)$, $A,B\leq G$ and $[A,B]=1$. Then $[A^{(2),\Omega},B^{(2),\Omega}]=1$. In particular, if $G$
is abelian then $G^{(2),\Omega}$ is also abelian. Furthermore,
if $H\leq G$ then $(C_G(H))^{(2),\Omega}\leq C_{G^{(2),\Omega}}(H^{(2),\Omega})$.
\end{lemma}
\proof
Let $x\in A^{(2),\Omega}$ and $y\in B^{(2),\Omega}$. Then for each $\alpha\in\Omega$ there exist $g_\alpha\in A$ and $h_\alpha\in B$ such that
$(\alpha,\alpha^x)^y=(\alpha,\alpha^x)^{g_\alpha}$ and $(\alpha,\alpha^y)^x=(\alpha,\alpha^y)^{h_\alpha}$ which imply
that $\alpha^{xy}=\alpha^{xg_\alpha}$, $\alpha^y=\alpha^{g_\alpha}$, $\alpha^x=\alpha^{h_\alpha}$ and
$\alpha^{yx}=\alpha^{yh_\alpha}$. So
$\alpha^{xy}=\alpha^{h_\alpha g_\alpha}=\alpha^{g_\alpha h_\alpha }=\alpha^{y h_\alpha}=\alpha^{yx}$. Thus
$xy=yx$ which means that $[A^{(2),\Omega},B^{(2),\Omega}]=1$.
\eqd

\begin{remark}\label{rmk}
{\rm Let $G\leq\Sym(\Omega)$. By Lemma \ref{commute},
$Z(G)\leq Z(G)^{(2),\Omega}\leq Z(G^{(2),\Omega})$. Furthermore,
if $G$ is $2$-closed on $\Omega$ then the centralizer
of each subgroup of $G$ is $2$-closed on $\Omega$. In particular, the center of $G$ is a $2$-closed group on  $\Omega$.}
\end{remark}

The proof of the following lemma is straightforward and so we omit the proof.
\begin{lemma}\label{permiso}
Let $G\leq\Sym(\Omega)$ and $H\leq\Sym(\Gamma)$ be permutation isomorphic. Then
$G^{(2),\Omega}\leq\Sym(\Omega)$ and $H^{(2),\Gamma}\leq\Sym(\Gamma)$ are permutation isomorphic.\eqd
\end{lemma}

\begin{lemma}\label{conj}
Let $G\leq\Sym(\Omega)$ and $x\in\Sym(\Omega)$. Then $(x^{-1}Gx)^{(2),\Omega}=x^{-1}G^{(2),\Omega}x$. In particular, $N_{\Sym(\Omega)}(G)\leq N_{\Sym(\Omega)}(G^{(2),\Omega})$.
\end{lemma}
\proof
Let $a\in (x^{-1}Gx)^{(2),\Omega}$ and $\alpha,\beta\in\Omega$. Then there exists $g\in G$ such that
$(\alpha^x,\beta^x)^a=(\alpha^x,\beta^x)^{x^{-1}gx}=(\alpha,\beta)^{gx}$. Thus $(\alpha,\beta)^{xax^{-1}}=(\alpha,\beta)^g$,
which implies that $xax^{-1}\in G^{(2),\Omega}$. Hence $a\in x^{-1}G^{(2),\Omega}x$ and so
$(x^{-1}Gx)^{(2),\Omega}\subseteq x^{-1}G^{(2),\Omega} x$.
Conversely, let $a\in x^{-1}G^{(2),\Omega}x$. Then $xax^{-1}=g$ for some $g\in G^{(2),\Omega}$.
Let $\alpha,\beta\in\Omega$. Then there exists
$h\in G$ such that $(\alpha^{x^{-1}},\beta^{x^{-1}})^g=(\alpha^{x^{-1}},\beta^{x^{-1}})^h$.
So $(\alpha,\beta)^a=(\alpha,\beta)^{x^{-1}hx}$, which means that $a\in (x^{-1}Gx)^{(2),\Omega}$. Hence $ x^{-1}G^{(2),\Omega}x\subseteq (x^{-1}Gx)^{(2)}$ and the proof is complete.
\eqd

\begin{corollary}\label{permiso2}
Let $G\leq\Sym(\Omega)$ and $H\leq\Sym(\Gamma)$ be permutation isomorphic. Then $G^{(2),\Omega}=G$ if and only if $H^{(2),\Gamma}=H$.
\end{corollary}
\proof Let $G^{(2),\Omega}=G$. Then by Lemma \ref{permiso}, $H$ and $H^{(2),\Gamma}$ are permutation isomorphic, as subgroups
of $\Sym(\Gamma)$. Now \cite[Exercise 1.6.1]{Dixon} implies that $H$ and $H^{(2),\Gamma}$ are conjugate in $\Sym(\Gamma)$.
Since, by \cite[Theorem 5.9]{Wielandt}, $(H^{(2),\Gamma})^{(2),\Gamma}=H^{(2),\Gamma}$, Lemma \ref{conj} implies that
$H^{(2),\Gamma}=H$. Similary, one can prove the converse direction.
\eqd

In \cite{math}, D. Holt pointed out that cyclic groups of prime power order and generalized quaternion groups are $2$-closed groups. We prove
the latter for completeness.
\begin{lemma}\label{holt}  Finite cyclic groups of prime power orders and generalized quaternion groups are $2$-closed.
\end{lemma}
\proof
Let $G\leq\Sym(\Omega)$ for some set $\Omega$. First, let $G$ be a finite cyclic $p$-group, $p$ a prime. We claim that $G_\alpha=1$ for some $\alpha\in\Omega$.
 Then \cite[Theorem 5.12]{Wielandt}
implies that $G$ is $2$-closed on $\Omega$.  Suppose, towards a contradiction, that for all $\alpha\in\Omega$, $G_\alpha\neq 1$. Thus for each
$\alpha\in\Omega$ there exists an element $x_\alpha\in G_\alpha$ of order $p$. Since $G$ has only one subgroup of order $p$, there
exists $x\in G$ of order $p$ such that for each $\alpha\in\Omega$, $\langle x\rangle=\langle x_\alpha\rangle$. Hence
$\langle x\rangle\leq \bigcap_{\alpha\in\Omega}G_\alpha(=1)$, a contradiction.

Since the generalized quaternion group $Q_{2^n}$, $n\geq 3$, has only one subgroup of order 2, the later case can be proved in
a similar argument.
\eqd

\begin{proposition}\label{direct}
Let $G=H\times K\leq\Sym(\Omega)$, where $(|H|,|K|)=1$.
Then $Z(G)^{(2),\Omega}=Z(H)^{(2),\Omega}\times Z(K)^{(2),\Omega}$.
\end{proposition}
\proof
Let $x\in Z(G)^{(2),\Omega}$ and $\alpha,\beta\in\Omega$. Then there exists $g=hk\in Z(G)$ such that $(\alpha,\beta)^x=(\alpha,\beta)^{hk}$, where $h\in Z(H)$ and $k\in Z(K)$.
Let $n=|H|$ and $m=|K|$. Then there exist $r,s\in\Bbb Z$
such that $rn+sm=1$. Since, by Lemma \ref{commute}, $Z(G)^{(2),\Omega}$ is abelian, $\alpha^{x^n}=\alpha^{h^nk^n}=\alpha^{k^n}$ and $\beta^{x^n}=\beta^{h^nk^n}=\beta^{k^n}$, which implies that $x^n\in Z(K)^{(2),\Omega}$.  Similarly, we have $x^m\in Z(H)^{(2),\Omega}$. Hence
$x=x^{rn+sm}=x^{rn}x^{sm}\in Z(H)^{(2),\Omega}Z(K)^{(2),\Omega}$.
Now let $x\in Z(H)^{(2),\Omega}\cap Z(K)^{(2),\Omega}$. Then for each
$\alpha\in\Omega$ there exist $h\in Z(H)$ and $k\in Z(K)$ such that $\alpha^x=\alpha^h=\alpha^k$. Thus $hk^{-1}\in Z(G)_\alpha$. Hence
by  Lemma \ref{stab}, $h\in Z(H)_\alpha$. Thus $\alpha^x=\alpha$. This proves $x=1$ and completes the proof.
\eqd


\begin{remark}{\rm Note that if $(|H|,|K|)\neq 1$, then the conclusion of Proposition \ref{direct} may not be true. For, let $G=\langle (1,2)(3,4),(3,4)(5,6)\rangle$.
 Then $G^{(2),\Omega}=\langle (1,2),(3,4),(5,6)\rangle$, where $\Omega=\{1,\ldots,6\}$. But $A^{(2),\Omega}=A$ and
 $B^{(2),\Omega}=B$, where $A=\langle (1,2)(3,4)\rangle$ and $B=\langle (3,4)(5,6)\rangle$.}
\end{remark}
The following corollary is proved in \cite[Theorem 12]{Monks}. Our proof is slightly different.
\begin{corollary}{\rm (\cite[Theorem 12]{Monks})}\label{cyclic}
Every finite cyclic group is a $2$-closed group.
\end{corollary}
\proof
Let $G$ be a finite cyclic permutation group of order $n=p_1^{n_1}\ldots p_k^{n_k}$, where $p_i$'s are distinct primes
and $n_i$'s are positive integers. Then $G\cong C_{p_1^{n_1}}\times\cdots\times C_{p_k^{n_k}}$. Now the result follows from
 Lemma \ref{holt} and Proposition \ref{direct} immediately.
\eqd

\begin{lemma}\label{abp}
Let $G$ be a finite abelian $p$-group, $p$ a prime. Then $G$ is a $2$-closed group if and only if it is cyclic.
\end{lemma}
\proof If $G$ is cyclic, then the result follows from Corollary \ref{cyclic}. Conversely, suppose that $G$ is a closed group. We prove
that $G$ is cyclic. Suppose, towards a contradiction, that $G$ is not cyclic. Then $G=G_1\times \cdots\times G_n$, where
$G_i\cong C_{p^{k_i}}$, $1\leq k_1\leq k_2\leq\cdots\leq k_n$ and $n\geq 2$. Let $G_i=\langle g_i\rangle$, $i=1,\ldots,n$. Then
$G=\langle g_1,\ldots,g_n\rangle$.  Let $\Omega_i=\{\alpha^i_{1},\ldots,\alpha^i_{p^{k_i}}\}$,
$i=1,\ldots,n$ and $\Omega_{n+1}=\{\beta_1,\ldots,\beta_p\}$ be $n+1$ pairwise disjoint sets and
 $\Omega=\bigcup_{i=1}^{n+1}\Omega_i$. Let $H=\langle h_1,\ldots,h_n\rangle$, where
\begin{eqnarray*}
h_1&=&(\beta_1,\ldots,\beta_p)(\alpha_1^1,\ldots,\alpha_{p^{k_1}}^1)\\
h_i&=&(\alpha_1^{i-1},\ldots,\alpha_{p^{k_{i-1}}}^{i-1})(\alpha_1^{i},\ldots,\alpha_{p^{k_{i}}}^{i}),~~~2\leq i\leq n.
\end{eqnarray*}

Then $H=\langle h_1\rangle\times\cdots\times\langle h_n\rangle\leq\Sym(\Omega)$. Since for each $i$, $g_i$ and $h_i$ has same order, there exists an isomorphism $\varphi:G\rightarrow H$. Hence
$G$ acts on $\Omega$ by the rule $\alpha^g=\alpha^{\varphi(g)}$. Furthermore, $G$ is also faithful on $\Omega$. Hence $G$
and $H$ are permutation isomorphic on $\Omega$. Thus Lemma \ref{permiso2} implies that $H^{(2),\Omega}=H$ if and only if
$G^{(2),\Omega}=G$. Since $G$ is a $2$-closed group, $G^{(2),\Omega}=G$ and so $H^{(2),\Omega}=H$. Let 
$x=(\beta_1,\ldots,\beta_p)$ and
$\gamma,\delta\in\Omega$. Then
 \begin{eqnarray*}
(\gamma,\delta)^{x}&=&\left\{
\begin{array}{ll}
(\gamma,\delta)^{h_1} & \gamma,\delta\in\Omega_{n+1}\\
(\gamma,\delta)& \gamma,\delta\notin\Omega_{n+1}\\
(\gamma,\delta)^{h_1h_2^{-1}}~ {\textrm or}~(\gamma,\delta)^{h_1}& \gamma\in\Omega_{n+1},\delta\notin\Omega_{n+1}\\
(\gamma,\delta)^{h_1h_2^{-1}}~ {\textrm or}~(\gamma,\delta)^{h_1}& \gamma\notin\Omega_{n+1}, \delta\in\Omega_{n+1}
\end{array}
\right.,
\end{eqnarray*}
which means that $x\in H^{(2),\Omega}$. Since $H\leq H^{(2),\Omega}$, for all $i=1,\ldots, n$, we have
$(\alpha_1^{i},\ldots,\alpha_{p^{k_{i}}}^{i})\in H^{(2),\Omega}$. Since $H$ is abelian, Lemma \ref{commute} implies that
$H^{(2),\Omega}$ is also abelian and so
$|H^{(2)}|\geq p.p^{k_1}p^{k_2}\cdots p^{k_n}=p|H|$, a contradiction.
\eqd

\begin{lemma}\label{disjoint} Let $G_i\leq\Sym(\Omega_i)$, $\Omega$ be disjoint union of $\Omega_i$'s, $i=1,\ldots,n$ and
$G=G_1\times\cdots\times G_n$. Then the natural action  of $G$ on $\Omega$ is faithful. Furthermore, if $G^{(2),\Omega}=G$
then for each $i=1,\ldots,n$,
$G_i^{(2),\Omega_i}=G_i$. In particular, if $G$ is a $2$-closed group then $G_i$ is a $2$-closed group, $i=1,\ldots,n$.
\end{lemma}
\proof The first part is clear.
Let $i\in\{1,\ldots,n\}$ and $x\in G_i^{(2),\Omega_i}$. Then for each $\alpha,\beta\in\Omega_i$ there exists
$g_{\alpha,\beta}\in G_i$ such that
$(\alpha,\beta)^x=(\alpha,\beta)^{g_\alpha,\beta}$. Put $\overline{x}=(x_1,\ldots,x_n)$, where $x_i=x$
 and for all $j\neq i$, $x_j=1$. Let $\gamma,\delta\in\Omega$. Then
 \begin{eqnarray*}
(\gamma,\delta)^{\overline{x}}&=&\left\{
\begin{array}{ll}
(\gamma,\delta)^{x} & \gamma,\delta\in\Omega_i\\
(\gamma,\delta)& \gamma,\delta\notin\Omega_i\\
(\gamma^x,\delta)& \gamma\in\Omega_i,\delta\notin\Omega_i\\
(\gamma,\delta^x)& \gamma\notin\Omega_i, \delta\in\Omega_i
\end{array}
\right.\\
&=&(\gamma,\delta)^{\overline{g}},
\end{eqnarray*}
where $\overline{g}=(g_1,\ldots,g_n)$, $g_i=g_{\gamma,\delta}$ and for all $j\neq i$, $g_j=1$. This implies that
$\overline{x}\in G^{(2),\Omega}$. Hence, $\overline{x}\in G$ and $x\in G_i$. Thus $G_i^{(2),\Omega_i}=G_i$.
\eqd \\

We now show that finite cyclic groups are the only finite abelian $2$-closed groups.

\begin{theorem}\label{ab} A finite abelian  group is $2$-closed if and only if it is cyclic.
\end{theorem}
\proof
By Corollary \ref{cyclic}, it is enough to prove that if $G$ is not cyclic, then  it is not a $2$-closed group. Suppose, towards a contradiction, that
$G$ is a non-cyclic $2$-closed group. Let
$|G|=p_1^{n_1}\cdots p_k^{n_k}$, where $p_i$'s are distinct primes and $n_i\geq 1$. Then $G=P_1\times\cdots\times P_k$,
where $P_i$ is the Sylow $p_i$-subgroup of $G$.
Since $G$ is not cyclic, at least one of $P_i$'s, say $P_1$, is not cyclic. Hence, by Lemma \ref{abp}, there exists a set $\Omega_1$
such that $P_1\leq\Sym(\Omega_1)$ and $P_1^{(2),\Omega_1}\neq P_1$, which contradicts Lemma \ref{disjoint}.
\eqd

\subsection{Proof of Theorems 1 and 2}

In this section, we are going to classify all finite nilpotent $2$-closed groups. Let us recall an important result:

\begin{theorem}{\rm (Universal embedding theorem \cite[Theorem 2.6 A]{Dixon})}\label{universal} Let $G$ be an arbitrary group
with a normal subgroup $N$ and put $K:=G/N$. Let $\psi: G\rightarrow K$ be a homomorphism of $G$ onto $K$ with kernel $N$.
Let $T:=\{t_u\mid u\in K\}$ be a set of right coset representatives of $N$ in $G$ such that $\psi(t_u)=u$ for each $u\in K$. Let
$x\in G$ and $f_x:K\rightarrow N$ be the map with $f_x(u)=t_uxt_{u\psi(x)}^{-1}$ for all $u\in K$. Then $\varphi(x):=(f_x,\psi(x))$
defines an embedding $\varphi$ of $G$ into $N\wr K$. Furthermore, if $N$ acts faithfully on a set $\Delta$ then $G$ acts
faithfully on $\Delta\times K$ by the rule $(\delta,k)^{x}=(\delta^{f_x(k)},k\psi(x))$.
\end{theorem}

Now we are ready to prove Theorem \ref{cycliccenter}.

{\bf Proof of Theorem \ref{cycliccenter}.} Suppose, towards a contradiction, that $Z(G)$ is not cyclic. Then there exists a non-cyclic 
Sylow $p$-subgroup $N$ of $Z(G)$. Hence, by Lemma \ref{abp}, there exists a set $\Delta$ such that 
$N\leq\Sym(\Delta)$ and $N^{(2),\Delta}\neq N$.
Let $K=G/N$ and  $\Gamma=\Delta\times K$. Keeping the notations of Theorem \ref{universal}, $G$ acts faithfully on $\Gamma$ by the rule
\[(\delta,k)^x=(\delta^{f_x(k)},k\psi(x))~~~\forall\delta\in\Delta,\forall k\in K,\forall x\in G.\]
Furthermore, since $N$ is in the center of $G$,
if $a\in N$ then for all $\delta\in\Delta$ and $k\in K$, we have $(\delta,k)^a=(\delta^a,k)$.



 Since $G$ is a
 $2$-closed group, Remark \ref{rmk} implies that $Z(G)$ is $2$-closed on $\Gamma$. Hence $N$ is $2$-closed on $\Gamma$,
 by \cite[Exercise 5.28]{Wielandt}. Let $\theta\in N^{(2),\Delta}$. Then for each $\delta_1,\delta_2\in\Delta$ there
 exists $x_{\delta_1,\delta_2}\in N$ such that
 $(\delta_1,\delta_2)^\theta=(\delta_1,\delta_2)^{x_{\delta_1,\delta_2}}$. Define $\overline{\theta}:\Gamma\rightarrow\Gamma$
 by the rule
  $(\delta,k)^{\overline{\theta}}=(\delta^\theta,k)$. Then $\overline{\theta}\in\Sym(\Gamma)$ and for each
  $(\delta_1,k_1),(\delta_2,k_2)\in\Gamma$, we have
  \[((\delta_1,k_1),(\delta_2,k_2))^{\overline{\theta}}=((\delta_1^\theta,k_1),(\delta_2^\theta,k_2))
  =((\delta_1^{x_{\delta_1,\delta_2}},k_1),(\delta_2^{x_{\delta_1,\delta_2}},k_2))
  =((\delta_1,k_1),(\delta_2,k_2))^{x_{\delta_1,\delta_2}},\]
which means that $\overline{\theta}\in N^{(2),\Gamma}$. Hence $\overline{\theta}\in N$ and so
$\overline{\theta}=y$ for some $y\in N$. Thus for each $\delta\in\Delta$ and $k\in K$, $(\delta,k)^{\overline{\theta}}=(\delta,k)^y$, which
implies that for each $\delta\in\Delta$, $\delta^\theta=\delta^y$. Hence $\theta=y\in N$, which proves that $N$ is $2$-closed on
$\Delta$, a contradiction.\eqd

Now we study finite $2$-closed $2$-groups. First, we consider finite $2$-groups having a normal elementary 
abelian subgroup of rank $2$.

\begin{lemma}\label{2-group} Let $G$ be a finite $2$-group having a normal subgroup $N\cong C_2\times C_2$. Then $G$ is not a $2$-closed
group.
\end{lemma}
\proof Suppose, towards a contradiction, that $G$ is a $2$-closed group. Without loss of generality, we may assume that
$N=\langle (1,2),(3,4)\rangle$. Let $a=(1,2)$ and $b=(3,4)$. Since, by Theorem 2, $Z(G)$ is cyclic, we may assume that 
$a\in Z(G)$ and we have $b\notin Z(G)$. 

Let $C=C_G(N)$ and $L=C/N$. Then $|G:C|=2$ and $G=C\langle t\rangle$ for some $t\in G\setminus C$, where $t^2\in C$.
Since $b\notin Z(G)$, $t^{-1}bt=ab$.
Let $\Delta=\{1,2,3,4\}$. Then, by Universal embedding theorem \ref{universal},
$C$ embeds in $N\wr C/N \leq \Sym(\Delta \times C/N)$ . Let $\overline{s}\in C/N$. Then for each $i\in\{1,\dots,4\}$ and $x\in N$, 
$(i,\overline{s})^x=(i^x,\overline{s})$ and moreover
\[ C_{(1,\overline{s})}=C_{(2,\overline{s})}=\langle b\rangle,~~C_{(3,\overline{s})}=C_{(4,\overline{s})}=\langle a\rangle.\]

Now, again by Universal embedding theorem, $G$ embeds in $C\wr G/C \leq \Sym(\Gamma \times G/C)$, where $\Gamma:=\Delta\times C/N$. Let $\{\overline{1},\overline{t}\}$ be the set of right transversal
of $C$ in $G$. Then for each $\gamma\in \Gamma$, we have 
\begin{eqnarray*}
 &&(\gamma,\overline{1})^g=\left\{
\begin{array}{ll}
(\gamma^g,\overline{1}) &g\in C\\
(\gamma^{gt^{-1}},\overline{t})&g\notin C
\end{array}
\right.,\\
&&(\gamma,\overline{t})^g=\left\{
\begin{array}{ll}
(\gamma^{tgt^{-1}},\overline{t}) &g\in C\\
(\gamma^{tg},\overline{1})&g\notin C
\end{array}
\right..
 \end{eqnarray*}
 Then for each $\overline{s}\in C/N$, since $t^{-1}bt=ab$, it follows that
 \[G_{((1,\overline{s}),\overline{1})}=G_{((2,\overline{s}),\overline{1})}=\langle b\rangle,~~
 G_{((1,\overline{s}),\overline{t})}=G_{((2,\overline{s}),\overline{t})}=\langle ab\rangle, \]
 and
 \[G_{((3,\overline{s}),\overline{1})}=G_{((3,\overline{s}),\overline{t})}=G_{((4,\overline{s}),\overline{1})}=
 G_{((4,\overline{s}),\overline{t})}=\langle a\rangle.\]
 
 Let $\Omega=\Gamma\times G/C$, $\overline{s}\in C/N$ and $\overline{k}\in G/C$. Define the 
 map $\theta$ on $\Omega$  as follows
 \begin{eqnarray*}
 ((i,\overline{s}),\overline{k})\mapsto \left\{
\begin{array}{ll}
((i,\overline{s}),\overline{k}) &i=1,2\\
((4,\overline{s}),\overline{k}) &i=3\\
((3,\overline{s}),\overline{k}) & i=4
\end{array}
\right..
 \end{eqnarray*}
 
 Then clearly $1\neq \theta\in\Sym(\Omega)$. We claim that $\theta\in G^{(2),\Omega}$. Let $\alpha=((i,\overline{s}),\overline{k})$
 and $\beta=((j,\overline{s'}),\overline{k'})$. Then
 \begin{eqnarray*}
 (\alpha,\beta)^\theta=\left\{
\begin{array}{ll}
(\alpha,\beta) & i,j\in\{1,2\}\\
(\alpha,\beta)^b & i,j\in\{3,4\}\\
(\alpha,\beta)^{b}~{\rm or}~(\alpha,\beta)^{ab} & i\in\{1,2\},j\in\{3,4\}\\
(\alpha,\beta)^{b}~{\rm or}~(\alpha,\beta)^{ab}  & i\in\{3,4\},j\in\{1,2\}
\end{array}
\right..
 \end{eqnarray*}
 This proves our claim. Now $\theta$ fixes $((1,\overline{1}),\overline{1})$ and $((1,\overline{1}),\overline{t})$.
 Since $G_{((1,\overline{1}),\overline{1})}\cap G_{((1,\overline{1}),\overline{t})}=\langle b\rangle \cap\langle ab\rangle=1$,
 we conclude that $G^{(2),\Omega}\neq G$, a contradiction. 
 \eqd \\
 
 In the following lemma, we prove that no finite nilpotent $2$-closed group can be splitted over any normal subgroup with an abelian complement.  
 
 \begin{lemma}\label{semidirectab} Let $G$ be a finite nilpotent group, $G=MH$, $M\unlhd G$, $H\leq G$ be abelian and $H\cap M=1$.
 Then $G$ is not a $2$-closed group.
\end{lemma}
\proof Suppose, towards a contradiction, that $G$ is a $2$-closed group. Then, by Theorem \ref{cycliccenter}, $G$ has a unique
minimal normal subgroup. Hence $H_G=1$.
Let $\Omega=\{Hm\mid m\in M\}$ be the set of right cosets of $H$ in $G$. Then 
$G$ acts on $\Omega$, by right multiplication, faithfully. Let
$\varphi:G\rightarrow\Sym(\Omega)$ be the corresponding monomorphism.

Let $H=\langle h_1\rangle\times\langle h_2\rangle\times\cdots\times\langle h_t\rangle$, $t\geq 1$ and $o(h_i)=n_i$. Let
$\Gamma_i$ be a set of $n_i$ elements such that $\Gamma_i\cap \Omega=\emptyset$, $\Gamma_i$'s be pairwise disjoint, 
$\tau_i$ be an $n_i$-cycle on $\Gamma_i$,
$\Omega'=(\bigcup_{i=1}^t\Gamma_i)\cup\Omega$ and 
$\overline{G}=\langle \varphi(M),\tau_1\varphi(h_1),\ldots, \tau_t\varphi(h_t)\rangle$.
Then $\overline{G}\leq\Sym(\Omega')$. 

Let $\{m_1,\ldots,m_l\}$ be a generating set of $M$. Then $G=\langle m_1,\ldots,m_l,h_1,\ldots,h_t\rangle$, where 
\begin{eqnarray*}
 h_i^{-1}m_jh_i=\omega_{i,j}(m_1,\ldots,m_l),
 \end{eqnarray*}
 $\omega_{i,j}(m_1,\ldots,m_l)$ is a word on $\{m_1,\ldots,m_l\}$, and we have some relations 
$r_m(m_1,\ldots,m_l)$, and $R_n(h_1,\ldots,h_t)~m=1,\ldots,s, n=1,\ldots,s'$, in terms of $m_1,\ldots,m_l$ and $h_1,\ldots,h_t$,
respectively. So 
\[\overline{G}=\langle \varphi(m_1),\ldots,\varphi(m_l),\tau_1\varphi(h_1),\ldots,\tau_t\varphi(h_t)\rangle\] 
and
\begin{eqnarray*}
(\tau_i\varphi(h_i))^{n_i}=1,(\tau_i\varphi(h_i))^{-1}\varphi(m_j)\tau_i\varphi(h_i)=\omega'_{i,j}(\varphi(m_1),\ldots,\varphi(m_l)),\\r'_m(\varphi(m_1),\ldots,\varphi(m_l)),R'_n(\varphi(h_1),\ldots,\varphi(h_t))
\end{eqnarray*}
where $\omega'_{i,j}(\varphi(m_1),\ldots,\varphi(m_l))=\varphi(\omega_{i,j}(m_1,\ldots,m_l))$,
$r'_m(\varphi(m_1),\ldots,\varphi(m_l))=\varphi(r_m(m_1,\ldots,m_l))$ and $R'_n(\varphi(h_1),\ldots,\varphi(h_t))=\varphi(R_n(h_1,\ldots,h_t))$. Hence, by Von Dyck's Theorem, there exists an epimorphism from
$G$ to $\overline{G}$. Let $\overline{H}=\langle \tau_1\varphi(h_1),\ldots,\tau_t\varphi(h_t)\rangle$. Then it is easy to see that 
$\overline{G}=\varphi(M)\overline{H}$ and 
$\varphi(M)\cap\overline{H}=1$. Hence $|\overline{G}|=|\varphi(M)||H|=|M||H|=|G|$, which
implies that $\overline{G}\cong G$. So
 $G$ acts faithfully on $\Omega'$.
 Clearly $\tau_1\notin\overline{G}$. We claim that $\tau_1$ is an element of the $2$-closure of $\overline{G}$ on $\Omega'$. Let 
 $\alpha,\beta\in\Omega'$. Then 
 \begin{eqnarray*}
(\alpha,\beta)^{\tau_1}&=&\left\{
\begin{array}{ll}
(\alpha,\beta) & \alpha,\beta\in\Omega\\
(\alpha,\beta)~{\textrm or} ~(\alpha,\beta)^{\tau_1\varphi(h_1)} & \alpha,\beta\notin\Omega\\
(\alpha,\beta)^{\varphi(m^{-1}h_1mh_1^{-1})\tau_1\varphi(h_1)}& \alpha=Hm\in\Omega,\beta\notin\Omega\\
(\alpha,\beta)^{\varphi(m^{-1}h_1mh_1^{-1})\tau_1\varphi(h_1)}& \alpha\notin\Omega,\beta=Hm\in\Omega,\\
\end{array}
\right.,
\end{eqnarray*}
which proves the claim. Hence $\overline{G}$ is not $2$-closed on $\Omega'$, a contradiction.
\eqd

\begin{corollary}\label{2-2} Let $G$ be a finite $2$-group. 
Then $G$ is a $2$-closed group if and only if it is cyclic or a generalized quaternion group.
\end{corollary}
\proof Let $G$ be a $2$-closed group. Then, by Lemma \ref{2-group}, $G$ has no  normal subgroup isomorphic to $C_2\times C_2$.
Hence, by \cite[Lemma 1.4]{Berkovich}, $G$ is cyclic, dihedral, semidihedral or a generalized quaternion group.
By Lemma \ref{semidirectab}, $G$ is not a dihedral or semidihedral group. Thus $G$ is cyclic or a generalized quaternion group.
The converse  follows from Lemma \ref{holt} and Corollary \ref{cyclic}.
\eqd

Now we study finite $2$-closed $p$-groups, where $p$ is an odd prime.

\begin{lemma}\label{podd}
Let $p$ be an odd prime and $G$ be a finite $p$-group. Then $G$ is a $2$-closed group if and only if it is cyclic.
\end{lemma}
\proof One direction is clear by Corollary \ref{cyclic}. Now we prove the converse direction.
Suppose, towards a contradiction, that $G$ is a non-cyclic $2$-closed group. Since $p\neq 2$, \cite[Lemma 1.4]{Berkovich}
 implies that $G$
has a normal subgroup $N$ isomorphic to $C_p\times C_p$.
Let $N=\langle a\rangle\times\langle b\rangle$.
By Proposition \ref{cycliccenter}, $Z(G)$ is cyclic. So we may assume that $a\in Z(G)$ and $b\notin Z(G)$. Hence $N\cap Z(G)=\langle a\rangle$.
Put $C:=C_G(N)$. Then $|G:C|=p$ and there exists $t\in G\setminus C$ such that $G=\langle t\rangle C$ and $t^p\in C$.

Let $H:=\langle b\rangle$ and $\Omega=\{Hg\mid g\in G\}$. Then $H_G=1$ and
$\Omega=\bigcup_{i=0}^{p-1}\Omega_i$, where $\Omega_i=\{Ht^ix\mid x\in C\}$. Note that $\Omega_i\cap\Omega_j=\emptyset$
whenever $i\neq j$. We claim that in the action of $G$ on $\Omega$, by right multiplication, $G^{(2),\Omega}\neq G$. Since $|N|=p^2$,
$N\leq Z_2(G)$. Let $i\in\{0,\ldots,p-1\}$. Then $[t^{i-2},b^{-1}]\in N\cap Z(G)$. So
$[t^{i-2},b^{-1}]=a^{s_i}$ for some $s_i\in\{0,\ldots,p-1\}$. On the other hand, if $i\neq 2$ then $s_i\neq 0$. To see this,
let $i\neq 2$ and $s_i=0$. Then $t^{i-2}\in C_G(b)$.
Since $C\leq C_G(b)\leq G$ and $b\notin Z(G)$, we have $C=C_G(b)$. Thus $t^{i-2}\in C$. Now $(i-2,p)=1$ and $t^p\in C$
imply that $t\in C$, which is a contradiction. Hence for each $i\neq 2$, $(s_i,p)=1$ and there exist $k_i,l_i\in\Bbb Z$ such
 that $k_is_i+l_ip=1$, which implies
that $[t^{i-2},b^{-k_i}]=a$.

We define the map $\theta :\Omega\rightarrow\Omega$ as follows:

$$ 
 (Ht^ix)^\theta= \begin{cases}
 Ht^ix & \text{if} \; i\neq 2\\
 Ht^ixa & \text{if} \; i=2
\end{cases}
$$

where $x\in C, i=0,\ldots , p-1.$
It is clear that $1\neq \theta\in \Sym(\Omega)$.  Let $\alpha,\beta\in\Omega$,
where $\alpha=Ht^ix\in\Omega_i$ and $\beta=Ht^jy\in\Omega_j$.
 Then
\begin{eqnarray*}
(\alpha,\beta)^\theta=\left\{
\begin{array}{ll}
(\alpha,\beta) & i,j\neq 2\\
(\alpha,\beta)^a & i,j=2\\
(\alpha,\beta)^{t^{-j}b^{k_j}t^j} & i=2,j\neq 2\\
(\alpha,\beta)^{t^{-i}b^{k_i}t^i} & i\neq 2,j=2
\end{array}
\right..
\end{eqnarray*}
Note that in the last two equalities, we use the facts $x,y\in C=C_G(N)$, $t^{-r}b^{k_r}t^r\in N$  and $[t^{r-2},b^{-k_r}]=a$,
where $r\in\{i,j\}$.
So we have proved that $\theta\in G^{(2),\Omega}$. Furthermore, $\theta\in (G^{(2),\Omega})_H\cap (G^{(2),\Omega})_{Ht}$
and $G_{H}\cap G_{Ht}=H\cap t^{-1}Ht=1$ which implies that $G\neq G^{(2),\Omega}$,
a contradiction.
\eqd

{\bf Proof of ``if" part of Theorem \ref{nilp}.} 
 Let $G$ be a finite nilpotent $2$-closed group. Then Lemma \ref{disjoint}
implies that all Sylow subgroups of $G$ are $2$-closed. Now the result follows from Corollary \ref{2-2} and Lemma \ref{podd}.
\eqd \\

Now to prove Theorem \ref{nilp}, it is enough to prove that $Q_{2^n}\times C_m$, where $n\geq 3$ and $m\geq 3$ is an odd integer
is a $2$-closed group. Holt \cite{math} mentioned that he thinks the group $Q_{2^n}\times C_m$ is $2$-closed for all odd integers $m>0$. In what follows we give a proof of the latter.  Let us start with the following lemma.



\begin{lemma}\label{stab} Let $G=H\times K\leq\Sym(\Omega)$. If $(|H|,|K|)=1$, then for each $\alpha\in\Omega$,
$G_\alpha=H_\alpha\times K_\alpha$.
\end{lemma}
\proof Let $\alpha\in\Omega$. It is clear that $H_\alpha,K_\alpha\unlhd G_\alpha$ and $H_\alpha\cap K_\alpha=1$. So it is enough to
prove that $G_\alpha=H_\alpha K_\alpha$. Clearly $H_\alpha K_\alpha\subseteq G_\alpha$. Now let $x\in G_\alpha$. Then
there exist $h\in H$ and $k\in K$ such that $x=hk$ and $\alpha^x=\alpha$. Hence $\alpha^h=\alpha^{k^{-1}}$.
On the other hand, $hk^{-1}=k^{-1}h$. Thus $\alpha^{h^t}=\alpha^{k^{-t}}$ for each integer $t\geq 1$. Let $n=|H|$
and $m=|K|$. Then there exist $r,s\in\Bbb Z$ such that $rn+sm=1$. Hence
\[ \alpha^h=\alpha^{h^{rn+sm}}=\alpha^{h^{sm}}=\alpha,~~\alpha^k=\alpha^{k^{rn+sm}}=\alpha^{k^{rn}}=\alpha.\]
So $x\in H_\alpha K_\alpha$, which completes the proof.\eqd

\begin{lemma}\label{abd} Let $G=H\times K\leq\Sym(\Omega)$ be transitive and $\Omega=\alpha^G$. If $(|H|,|K|)=1$, then
the action of $G$ on $\Omega$ is
equivalent to the action of $G$ on $\Omega_1\times\Omega_2$, where $\Omega_1=\alpha^H$, $\Omega_2=\alpha^K$
and $G$ acts on $\Omega_1\times\Omega_2$ by the rule $(\alpha^h,\alpha^k)^{g}=(\alpha^{hh_1},\alpha^{kk_1})$, where $g=h_1k_1$.
\end{lemma}
\proof First note that, by Lemma \ref{stab}, $G_\alpha=H_\alpha\times K_\alpha$. This implies that the map
$\lambda:\Omega\rightarrow\Omega_1\times\Omega_2$, where $\lambda(\alpha^{hk})=(\alpha^h,\alpha^k)$ is well-defined and
one-to-one. Clearly $\lambda$ is onto. Hence $\lambda$ is a bijection. Now let $\beta\in\Omega$ and $x=hk\in G$,
where $h\in H$ and $k\in K$. Then there exists $g=h_1k_1\in G$, where $h_1\in H$ and $k_1\in K$, such that
$\beta=\alpha^g$. Furthermore,
\[\lambda(\beta^x)=\lambda(\alpha^{gx})=\lambda(\alpha^{h_1hk_1k})=(\alpha^{h_1h},\alpha^{k_1k})=(\alpha^{h_1},\alpha^{k_1})^{x}=\lambda(\beta)^x.\]
This completes the proof.
\eqd

Recall that $G,H\leq\Sym(\Omega)$ are called {\it $k$-equivalent} if $G$ and $H$ has same orbits on $\Omega\times\cdots\times\Omega$
 ($k$-time),
 see \cite[Definition 4.1]{Wielandt}.

\begin{lemma}{\rm (\cite[Lemma 2.3.2]{Xu})} \label{equivalent}
Let $G,H\leq\Sym(\Omega)$ be $k$-equivalent, $k\geq 1$. Let $N$ be an intransitive normal subgroup of both $G$ and $H$.
Let $\overline{\Omega}$ be the set of orbits of $N$ on $\Omega$. Then $\overline{G}=G/N$ and $\overline{H}=H/N$ are $k$-equivalent 
on $\overline{\Omega}$.
\end{lemma}

\begin{lemma}\label{qt}
Let $G=H\times K\leq\Sym(\Omega)$ be finite, $H,K\neq 1$ and $H$ be abelian.
Then $H\unlhd G^{(2),\Omega}$. Furthermore, if $(|H|,|K|)=1$ then $H$ acts intransitively on $\Omega$ and moreover
$G/H$ and $G^{(2),\Omega}/H^{(2),\Omega}$ both act faithfully on $\overline{\Omega}:=\{\alpha^H\mid\alpha\in\Omega\}$.
\end{lemma}
\proof Since $H\leq Z(G)$ and, by Remark \ref{rmk}, $Z(G)\leq Z(G^{(2),\Omega})$, we have $H\unlhd G^{(2),\Omega}$. Now
let $(|H|,|K|)=1$. Since $K\neq 1$, Lemma \ref{stab} implies that $H$ is intransitive on $\Omega$.
On the other hand, $G^{(2),\Omega}$ (and also $G$) acts  on $\overline{\Omega}$ by the rule $(\alpha^H)^x=(\alpha^x)^H$,
where $\alpha\in\Omega$ and
$x\in G^{(2),\Omega}$. To complete the proof, it is enough to prove that $H^{(2),\Omega}$ and $H$ are the kernels of
the actions of $G^{(2),\Omega}$
and $G$ on $\overline{\Omega}$, respectively.

Let $L$ be the kernel of the action of $G$ on $\overline{\Omega}$. Then $H\leq L$. Let $x\in L$ and $\alpha\in\Omega$ be arbitrary.
Then $(\alpha^x)^H=\alpha^H$ which implies that $x\in G_\alpha H$. On the other hand, by Lemma \ref{stab},
$G_\alpha=H_\alpha\times K_\alpha$.
So $x\in HK_\alpha$. This shows that for each $\alpha\in\Omega$, there exist $k_\alpha\in K_\alpha$ and $h_\alpha\in H$
such that $x=h_\alpha k_\alpha$. Since every element in $G$ has a unique presentation as a product of an
element of $H$ by an element of $K$, for each $\alpha,\alpha'\in\Omega$, we have
$k_\alpha=k_{\alpha'}$. Thus $k=k_\alpha$, for each $\alpha\in\Omega$, fixes all elements of $\Omega$. Hence $k=1$
which implies that $x\in H$. This proves that $H=L$.

Finally, let $M$ be the kernel of $G^{(2),\Omega}$ on $\overline{\Omega}$. We prove that $M=H^{(2),\Omega}$. Let
$x\in H^{(2),\Omega}$.
Then for each $\alpha\in\Omega$, there exists $h_\alpha\in H$ such that $\alpha^x=\alpha^h$. Thus
for each $\alpha\in\Omega$, $(\alpha^H)^x=(\alpha^x)^H=(\alpha^{h_\alpha})^H=\alpha^H$. Hence $H^{(2),\Omega}\leq M$.

Let $x\in M$ and $\alpha,\beta\in\Omega$. Then $x\in G^{(2),\Omega}$, $(\alpha^x)^H=\alpha^H$ and $(\beta^x)^H=\beta^H$. Hence
there exist $h_1,h_2\in H$ such that $\alpha^x=\alpha^{h_1}$ and $\beta^x=\beta^{h_2}$. Furthermore, there exists $g\in G$ such that
$\alpha^x=\alpha^g$ and $\beta^x=\beta^g$. Thus $gh_1^{-1}\in G_\alpha$ and $gh_2^{-1}\in G_\beta$. Since
$G_\alpha=H_\alpha K_\alpha$ and $G_\beta=H_\beta K_\beta$, there exist $r_1\in H_\alpha$, $s_1\in K_\alpha$ and $r_2\in H_\beta$
and $s_2\in K_\beta$ such that $gh_1^{-1}=r_1s_1$ and $gh_2^{-1}=r_2s_2$. Hence $r_1h_1s_1=r_2h_2s_2$. Again, since
every element of $G$ has a unique presentation as a product of an element of $H$ by an element of $K$, $s_1=s_2$ and $r_1h_1=r_2h_2$. Hence
\[\alpha^x=\alpha^{h_1}=\alpha^{r_1^{-1} r_2h_2}=\alpha^{r_2h_2},~
\beta^x=\beta^{h_2}=\beta^{r_2^{-1}r_1h_1}=\beta^{r_1h_1}=\beta^{r_2h_2},\]
which implies that $x\in H^{(2),\Omega}$. This completes the proof.
\eqd

\begin{proposition}\label{genholt}
Let $G=H\times K\leq\Sym(\Omega)$ be finite, $(|H|,|K|)=1$, $H$ be abelian, $H^{(2),\Omega}=H$ and $K^{(2),\overline{\Omega}}=K$,
where $\overline{\Omega}=\{\alpha^H\mid \alpha\in\Omega\}$.
Then $G^{(2),\Omega}=G$.
\end{proposition}
\proof
By Lemma \ref{qt}, $H^{(2),\Omega}=H$ implies that $G/H$ and $G^{(2),\Omega}/H$ both act faithfully on $\overline{\Omega}$.
Also Lemma \ref{equivalent} implies that $G^{(2),\Omega}/H\leq (G/H)^{(2),\overline{\Omega}}$. Since $G/H\cong K$ and $K^{(2),\overline{\Omega}}=K$, Lemma \ref{permiso} implies that $(G/H)^{(2),\overline{\Omega}}\cong K$. Also $K\cong G/H\leq G^{(2),\Omega}/H$. Thus
$|G^{(2),\Omega}/H|=|K|$, which implies that $|G^{(2),\Omega}|=|H||K|=|G|$. Hence $G^{(2),\Omega}=G$.
\eqd \\

We are now ready to complete the proof of Theorem \ref{nilp} by giving the proof of its ``only if" part. \\ 

{\bf Proof of ``only if" part of Theorem \ref{nilp}.} 
It follows from Proposition \ref{genholt} and Lemma \ref{holt}.
\eqd \\

\section*{Acknowledgments}
This research was financially supported by Iran National Science Foundation: INSF.
The first author was supported in part
by grant No. 95050219 from School of Mathematics, Institute for Research in Fundamental Sciences (IPM). The first author
was additionally financially supported
by the Center of Excellence for Mathematics at the University of Isfahan.

\end{document}